\documentclass[dspr]{apjrnl}
\usepackage{amssymb}
\usepackage{amsmath}
\usepackage{amsfonts}
\usepackage[doublespacing]{setspace}
\usepackage{cite}

{\theoremstyle{plain}
\newtheorem{theorem}{Theorem}

\newtheorem{lemma}{Lemma}
{\theoremstyle{remark}

}
{\theoremstyle{definition}

}

\iftwocol\AtEndDocument{\end{multicols}}\fi
\input{tcilatex}
\begin{document}

\title{Uniqueness of solution of an inverse source problem for
ultrahyperbolic equations}
\author{Fikret G\"{o}lgeleyen$^{\text{\ref{first}}}$ and Masahiro Yamamoto$^{%
\text{\ref{second}, \ref{third}}}$ \\
\additem[first]{Department of Mathematics, Faculty of Arts and Sciences, Zonguldak B\"ulent Ecevit University,
Zonguldak 67100 Turkey} \\
\additem[second]{Department of Mathematical Sciences,
The University of Tokyo,  3-8-1 Komaba, Meguro, Tokyo 153-8914 Japan} \\
\additem[third]{
Honorary Member of Academy of Romanian Scientists,\\
Splaiul Independentei Street, no 54,
050094 Bucharest Romania} \\
E-mail: f.golgeleyen@beun.edu.tr, myama@ms.u-tokyo.ac.jp}
\maketitle

\begin{abstract}
The aim of this article is to investigate the uniqueness of solution of an
inverse problem for ultrahyperbolic equations. We first reduce the inverse
problem to a Cauchy problem for an integro-differential equation and then by
using a pointwise Carleman type inequality we prove the uniqueness.
\end{abstract}

\section{Introduction and the main result}

In this article, we consider an inverse problem for an ultrahyperbolic
equation. One of our motivations to deal with this equation is its
interesting structure from the point of view of the theory of partial
differential equations. For instance, depending on the specific form of
initial conditions, solutions possess both hyperbolic and elliptic
properties (see \cite{Kos1}). Another motivation is recent discussions on
the possibility of physics in multiple time dimensions, (e.g., \cite{Ba, Sp,
Te}). Namely, in some superstring theories which attempt to unify the
general theory of relativity and the quantum mechanics, extra dimensions are
required for the consistency of theory. When the presence of more than one
temporal dimension is considered, the mathematical model occurs as an
ultrahyperbolic equation (e.g., \cite{Cr}). More precisely, the paper \cite%
{Cr} asserts that the equation in a form of 
\begin{equation*}
\partial _{xx}u(x,y_{1},...,y_{m})-\sum_{j=1}^{m}\partial
_{y_{j}y_{j}}u(x,y_{1},...,y_{m})=F(x,y_{1},...,y_{m})
\end{equation*}%
is of central physical importance, which describes the dynamical evolution
of many physical quantities of classical and quantum field theories
including the components of the electromagnetic fields in the case of a
single time dimension, while the equation in a form of 
\begin{equation*}
\sum_{i=1}^{n}\partial
_{x_{i}x_{i}}u(x_{1},...,x_{n},y_{1},...,y_{m})-\sum_{j=1}^{m}\partial
_{y_{j}y_{j}}u(x_{1},...,x_{n},y_{1},...,y_{m})=F(x_{1},...,x_{n},y_{1},...,y_{m})
\end{equation*}%
is fundamental where $x\in \mathbb{R}^{n}$ and $y\in \mathbb{R}^{m}$ are
respectively time-like variables and space-like variables. \newline

Let $n,m\geq 2$. Inspired by \cite{Cr, Sp, Te}, we here consider an
ultrahyperbolic equation in $u(x,y):=u(x_{1},...,x_{n},y_{1},...,y_{m})$,
which is associated with general geometry in the space variables $y$: 
\begin{align*}
& Lu(x,y)\equiv \sum_{i=1}^{n}\partial
_{x_{i}x_{i}}u(x,y)-%
\sum_{i,j=1}^{m}a_{ij}(x,y_{1},...,y_{m-1})u_{y_{i}y_{j}}(x,y) \\
+&
\sum_{i=1}^{n}a_{i}(x,y)u_{x_{i}}(x,y)+\sum%
\limits_{j=1}^{m}b_{j}(x,y)u_{y_{j}}(x,y)+a_{0}(x,y)u(x,y)
\end{align*}%
\begin{equation}
=f(x,y)g(x,y_{1},...,y_{m-1})  \tag{1}  \label{1}
\end{equation}%
in the domain $\Omega =D\times G$. Here $D\subset \mathbb{R}^{n}$ and $%
G\subset \mathbb{R}^{m}$ are bounded domains, and we assume that $\Omega
\subset \mathbb{R}^{n+m}$ is supported by the plane $x_{1}=0$, $G=G^{\prime
}\times I$ with an open interval $I$ and $G^{\prime }\subset \mathbb{R}%
^{m-1} $, and the coefficients are assumed to satisfy $a_{ij}\left(
x,y_{1},...,y_{m-1}\right) \in C^{2}(\overline{D\times G^{\prime }}),$ $%
a_{i}(x,y),b_{j}(x,y)\in C\left( \overline{\Omega }\right) $ $%
(i=0,1,...,n;j=1,...,m),$ $f\in C^{2}\left( \overline{\Omega }\right) $. 
\newline

The purpose of this article is to investigate the uniqueness of solution of
the following problem:

\textbf{Problem.}\newline
\textit{\ For given $u_{0}(x,y_{1},...,y_{m-1})$, find a pair of functions }$%
(u(x,y),g(x,y_{1},...,y_{m-1}))$\textit{\ in }$\Omega $\textit{\ satisfying
equation (1), Cauchy data}%
\begin{equation}
u\left( 0,x_{2},...,x_{n},y\right) =u_{x_{1}}\left(
0,x_{2},...,x_{n},y\right) =0  \tag{2}  \label{2}
\end{equation}%
\textit{and the additional information}%
\begin{equation}
u\left( x,y_{1},...,y_{m-1},0\right) =u_{0}\left( x,y_{1},...,y_{m-1}\right)
.  \tag{3}  \label{3}
\end{equation}%
This is an inverse problem of determining a factor $g$ which is independent
of the component $y_{m}$ of the source in (1) which causes the action under
consideration. This inverse problem is called an inverse source problem. 
\newline

Our main result is stated in Theorem 1:

\begin{theorem}
Let $f(x,y^{\prime },0)\neq 0$ and the functions $a_{ij}$ satisfy%
\begin{equation}
-\overset{m}{\underset{i,j=1}{\sum }}\frac{\partial a_{ij}}{\partial x_{1}}%
\xi ^{i}\xi ^{j}\geq \alpha _{1}|\xi |^{2},\alpha _{1}>0  \tag{4}  \label{4}
\end{equation}%
for any $\xi \in \mathbb{R}^{m}$ and $(x,y^{\prime })\in D\times G^{\prime }$%
. Then Problem has at most one solution $(u,g)$ such that $(u,g)\in
C^{2}\left( \Omega \right) \times C\left( D\times G^{\prime }\right) $.
\end{theorem}

Inverse problems for ultrahyperbolic equations were studied in \cite%
{Am1,Am2,BK, GolYama}, where the key method is based on weighted a priori
estimates and was firstly developed by Bukhgeim and Klibanov \cite{BK}. A
uniqueness theorem for ultrahyperbolic equations, is given by \cite{BK} for
a bounded domain with Dirichlet and Numann type condition on a part of the
boundary. In \cite{Am1} and \cite{Am2}, uniqueness is invesigated in an
unbounded domain with an additional information for the solution of direct
problem at $y=0.$ In \cite{GolYama}, H\"{o}lder stability estimates were
obtained in a bounded domain by some lateral boundary data. A major
difference of our work from the existing results is that, in Problem,
additional information is given only at $y_{m}=0.$

As for the direct problem (1) - (2) with given $fg$, it is known that the
problem of determination of the function $u$ from relations (1) and (2) is
ill-posed in the Hadamard sense (see \cite{LRS}, Chapter 4). By using the
mean-value theorem of Asgeirsson, it was shown by \cite{Cou} that the
existence of solutions fails if the initial conditions are not properly
prescribed. We refer to \cite{Bu, Di, Mur, Ow, Ro}, as for the uniqueness
results for various Cauchy, Dirichlet and Neumann problems for
ultrahyperbolic equations. Finally, in \cite{Cr} it is proved that under a
nonlocal constraint, the initial value problem is well-posed for initial
data given on a codimension-one hypersurface.

\section{Key Carleman estimate}

We set 
\begin{eqnarray*}
x &=&(x_{1},^{\prime }x)\in \mathbb{R}^{n},\quad ^{\prime
}x=(x_{2},...,x_{n})\in \mathbb{R}^{n-1}, \\
y &=&\left( y^{\prime },y_{m}\right) \in \mathbb{R}^{m}, \quad y^{\prime
}=\left( y_{1},...,y_{m-1}\right) \in \mathbb{R}^{m-1},
\end{eqnarray*}
and 
\begin{eqnarray*}
\partial _{x_{i}} &=&\frac{\partial }{\partial x_{i}},\text{ }\partial
_{y_{j}}=\frac{\partial }{\partial y_{j}},\text{ }\nabla _{^{\prime
}x}=\left( \partial _{x_{2}},\partial _{x_{3}},\cdots ,\partial
_{x_{n}}\right) , \\
\nabla _{y} &=&\left( \partial _{y_{1}},\partial _{y_{2}},\cdots ,\partial
_{y_{m}}\right) ,\text{ }\Delta _{x}=\sum\limits_{i=1}^{n} \partial_{x_ix_i},%
\text{ }\Delta _{^{\prime }x}=\sum\limits_{i=2}^{n} \partial_{x_ix_i}.
\end{eqnarray*}

In order to prove Theorem 1, the key tool is an Carleman type inequality
which will be presented in Lemma 1 below. First of all, we reduce equation
(1) to a more suitable form by introducing a new variable $\widetilde{x}_{1}=%
\sqrt{2x_{1}}-\eta _{0}$, that is, $x_{1}=\frac{1}{2}\left( \widetilde{x}%
_{1}+\eta _{0}\right) ^{2}$, where $2\eta _{0}=\min \left\{ \alpha
_{0},\gamma \right\} ,$ the parameters $\alpha _{0},$ $\gamma $ will be
specified later, $\eta _{0}>0.$ Without loss of generality, we assume that $%
\sqrt{2x_{1}}-\eta _{0}>0,$ i.e., $x_{1}>\frac{\eta _{0}}{2}$, and so we
have $\widetilde{x}_{1}>0$.

Then, for the new function $\widetilde{u}\left( \widetilde{x}_{1},^{\prime
}x,y\right) \equiv u\left( \frac{1}{2}\left( \widetilde{x}_{1}+\eta
_{0}\right) ^{2},^{\prime }x,y\right) $, by using the relations%
\begin{eqnarray*}
u_{x_{1}} &=&\widetilde{u}_{\widetilde{x}_{1}}\frac{d\widetilde{x}_{1}}{%
dx_{1}}=\widetilde{u}_{\widetilde{x}_{1}}\frac{1}{\widetilde{x}_{1}+\eta _{0}%
}; \\
u_{x_{1}x_{1}} &=&\widetilde{u}_{\widetilde{x}_{1}\widetilde{x}_{1}}\left( 
\widetilde{x}_{1}+\eta _{0}\right) ^{-2}-\widetilde{u}_{\widetilde{x}%
_{1}}\left( \widetilde{x}_{1}+\eta _{0}\right) ^{-3},
\end{eqnarray*}%
we have 
\begin{equation*}
\left( \widetilde{x}_{1}+\eta _{0}\right) ^{-2}\widetilde{u}_{\widetilde{x}%
_{1}\widetilde{x}_{1}}+\Delta _{^{\prime }x}\widetilde{u}-\overset{m}{%
\underset{i,j=1}{\sum }}\widetilde{a}_{ij}\widetilde{u}_{y_{i}y_{j}}+\overset%
{n}{\underset{i=2}{\sum }}\widetilde{a}_{i}\widetilde{u}_{x_{i}}+\overset{m}{%
\underset{j=1}{\sum }}\widetilde{b}_{j}\widetilde{u}_{y_{j}}+\widetilde{a}%
_{1}\widetilde{u}_{\widetilde{x}_{1}}+\widetilde{a}_{0}\widetilde{u}=%
\widetilde{f}\widetilde{g},
\end{equation*}%
where $\widetilde{a}_{ij}=a_{ij}(\frac{1}{2}\left( \widetilde{x}_{1}+\eta
_{0}\right) ^{2},^{\prime }x,y^{\prime }),$ $i,j=1,...,m;$ $\widetilde{a}%
_{i}=a_{i}(\frac{1}{2}\left( \widetilde{x}_{1}+\eta _{0}\right)
^{2},^{\prime }x,y),$ $i=0,2,3,...,n;$ $\widetilde{b}_{j}=b_{j}(\frac{1}{2}%
\left( \widetilde{x}_{1}+\eta _{0}\right) ^{2},^{\prime }x,y),$ $%
j=1,2,...,m;\ \widetilde{a}_{1}=a_{1}(\frac{1}{2}\left( \widetilde{x}%
_{1}+\eta _{0}\right) ^{2},^{\prime }x,y)-\left( \widetilde{x}_{1}+\eta
_{0}\right) ^{-3};$ $\widetilde{f}=f(\frac{1}{2}\left( \widetilde{x}%
_{1}+\eta _{0}\right) ^{2},^{\prime }x,y),$ $\widetilde{g}=g(\frac{1}{2}%
\left( \widetilde{x}_{1}+\eta _{0}\right) ^{2},^{\prime }x,y^{\prime }).$

For the sake of simplicity, let us denote $\widetilde{u},$ $\widetilde{a}%
_{ij},$ $\widetilde{x}_{1},$ $\widetilde{a}_{k},$ $\widetilde{b}_{s},$ $%
\widetilde{f}$, $\widetilde{g}$ by $u,$ $a_{ij},$ $x_{1},$ $a_{s},$ $b_{j}$, 
$f,$ $g$ respectively, where $i,j=2,...,m;$ $k=0,1,2,...,n;$ $s=1,2,...,m.$
Then we can write%
\begin{align}
& \left( x_{1}+\eta _{0}\right) ^{-1}u_{x_{1}x_{1}}+\left( x_{1}+\eta
_{0}\right) (\Delta _{^{\prime }x}u-\sum_{i,j=1}^{m}a_{ij}u_{y_{i}y_{j}}) 
\notag \\
& +\left( x_{1}+\eta _{0}\right)
(\sum_{i=1}^{n}a_{i}u_{x_{i}}+\sum_{j=1}^{m}b_{j}u_{y_{j}}+a_{0}u)  \notag \\
& =\left( x_{1}+\eta _{0}\right) fg.  \tag{5}
\end{align}

We set%
\begin{equation*}
L_{0}u\equiv \left( x_{1}+\eta _{0}\right) ^{-1}u_{x_{1}x_{1}}+(x_{1}+\eta
_{0})\left( \Delta _{^{\prime
}x}u-\sum_{i,j=1}^{m}a_{ij}u_{y_{i}y_{j}}\right) .
\end{equation*}%
We introduce%
\begin{equation*}
\Omega _{\gamma }=\left\{ (x,y);\text{ }x_{1}>0,\text{ }0<\delta x_{1}+\frac{%
1}{2}\sum_{i=2}^{n}(x_{i}-x_{i}^{0})^{2}+\frac{1}{2}%
\sum_{i=1}^{m-1}(y_{i}-y_{i}^{0})^{2}+\frac{1}{2}y_{m}^{2}<\gamma \right\} 
\text{,}
\end{equation*}%
where $0<\gamma <1,$ $\delta >4,$ $\left( x^{0},y^{0}\right) \in \overline{%
\Omega },$ $x^{0}=\left( 0,x_{2}^{0},...,x_{n}^{0}\right) ,$ $y^{0}=\left(
y_{1}^{0},y_{2}^{0},...,y_{m-1}^{0},0\right) .$

In $\Omega _{\gamma }$ we define the weight function 
\begin{equation*}
\chi =\exp \left( \lambda \psi ^{-\nu }\right) ,\text{ }\psi (x)=\delta
x_{1}+\frac{1}{2}\sum_{i=2}^{n}(x_{i}-x_{i}^{0})^{2}+\frac{1}{2}%
\sum_{i=1}^{m-1}(y_{i}-y_{i}^{0})^{2}+\frac{1}{2}y_{m}^{2}+\alpha _{0}\text{,%
}
\end{equation*}%
where $\alpha _{0}>0,$ $\gamma +\alpha _{0}=\rho <1,$ $\alpha _{0}<\psi
(x)<\rho ,$ and $\lambda ,$ $\nu ,$ $\delta $ are positive parameters
satisfying some additional conditions which are specified later.

The following Carleman estimate is the key for the proof of Theorem 1.

\begin{lemma}
Let condition (4) be satisfied. With arbitrarily fixed constant $M>0$, we
assume that%
\begin{equation*}
\left\vert \left\vert a_{ij}\right\vert \right\vert _{C^{2}(\overline{\Omega 
})}\leq M
\end{equation*}%
and the number $\gamma $ is "small", that is%
\begin{equation}
0<\gamma <\min \{\frac{1}{2},\frac{4}{3}(M\left( Mm^{2}\varepsilon
_{0}^{-1}+m^{2}+m(m+1)\right) )^{-1/2}\},  \tag{6}
\end{equation}%
where $0<\varepsilon _{0}<\frac{\alpha _{1}}{4m}$. Then there exists a
constant $\delta _{\ast }=\delta _{\ast }(\alpha _{1},M,n,m,\nu )>0$ such
that for any $\delta >\delta _{\ast }$ there exists $\lambda _{\ast
}=\lambda _{\ast }(\delta )$ such that the following estimate holds: 
\begin{eqnarray}
&&\psi ^{\nu +1}\left( L_{0}\varphi \right) ^{2}\chi ^{2}-2\lambda \nu \beta
_{0}\varphi \left( x_{1}+\eta _{0}\right) \left( L_{0}\varphi \right) \chi
^{2}  \notag \\
&\geq &2\lambda \nu \delta \left( x_{1}+\eta _{0}\right) ^{-3}\varphi
_{x_{1}}^{2}\chi ^{2}+2\lambda \nu \left( x_{1}+\eta _{0}\right)
^{2}\left\vert \nabla _{^{\prime }x}\varphi \right\vert ^{2}\chi ^{2}  \notag
\\
&&+2\lambda \nu \left( x_{1}+\eta _{0}\right) \left\vert \nabla _{y}\varphi
\right\vert ^{2}\chi ^{2}+\lambda ^{3}\nu ^{4}\delta ^{4}\psi ^{-2\nu
-3}\varphi ^{2}\chi ^{2}+D\left( \varphi \right)  \TCItag{7}
\end{eqnarray}%
for all $\varphi \left( x,y\right) \in C^{2}\left( \overline{\Omega }%
_{\gamma }\right) $ and $\lambda >\lambda _{\ast }$. In (7), $\beta
_{0}\equiv \beta _{0}(n,m)=n+2+Mm((1+3\sqrt{2\gamma })m+1)$ and $D\left(
\varphi \right) $ is described by a divergence form which includes the
function $\varphi $ and is given explicitly in the proofs of the lemmata
below.
\end{lemma}

\section{The proof of Lemma 1.}

In the proof of Lemma 1, we shall use two Lemmata 2 and 3.

\begin{lemma}
Under the hypothesis of Lemma 1, there exists a constant $\delta _{0}=\delta
_{0}(\alpha _{1},M,n,m,\nu )>0$ such that for any $\delta >\delta _{0}$
there exists $\lambda _{0}=\lambda _{0}(\delta )$ such that 
\begin{eqnarray}
\psi ^{\nu +1}\left( L_{0}\varphi \right) ^{2}\chi ^{2} &\geq &2\lambda \nu
\delta \left( x_{1}+\eta _{0}\right) ^{-3}\varphi _{x_{1}}^{2}\chi
^{2}-2\lambda \nu \left( x_{1}+\eta _{0}\right) ^{2}(\beta _{0}-1)\left\vert
\nabla _{^{\prime }x}\varphi \right\vert ^{2}\chi ^{2}  \notag \\
&&+\lambda \nu \delta \alpha _{1}\left( x_{1}+\eta _{0}\right) \left\vert
\nabla _{y}\varphi \right\vert ^{2}\chi ^{2}+2\delta ^{4}\lambda ^{3}\nu
^{4}\left( x_{1}+\eta _{0}\right) ^{-2}\psi ^{-2\nu -3}\varphi ^{2}\chi ^{2}
\notag \\
&&+\sigma _{1}(\lambda ,\delta )\varphi ^{2}\chi ^{2}+D_{1}(\chi \varphi
)+D_{2}(\chi \varphi ),  \TCItag{8}
\end{eqnarray}%
for all $\lambda >\lambda _{0}$ and $\varphi \in C^{2}(\overline{\Omega }%
_{\gamma })$. Here%
\begin{equation*}
\sigma _{1}(\lambda ,\delta )=\lambda ^{3}\nu ^{3}\sigma _{11}+\lambda
^{2}\nu ^{2}\sigma _{12},
\end{equation*}%
\begin{eqnarray*}
\sigma _{11} &=&\psi ^{-2\nu -2}(-2\delta ^{3}\left( x_{1}+\eta _{0}\right)
^{-3}+2\left( x_{1}+\eta _{0}\right) ^{2}(\beta _{0}-1)\left\vert \nabla
_{^{\prime }x}\psi \right\vert ^{2} \\
&&-\delta \alpha _{1}\left( x_{1}+\eta _{0}\right) \left\vert \nabla
_{y}\psi \right\vert ^{2}), \\
\sigma _{12} &=&(\nu +1)\psi ^{-\nu -2}(-2\delta ^{3}\left( x_{1}+\eta
_{0}\right) ^{-3}+2\left( x_{1}+\eta _{0}\right) ^{2}(\beta
_{0}-1)\left\vert \nabla _{^{\prime }x}\psi \right\vert ^{2} \\
&&-\delta \alpha _{1}\left( x_{1}+\eta _{0}\right) (\nu +1)\left\vert \nabla
_{y}\psi \right\vert ^{2})+\psi ^{-\nu -1}(6\delta ^{2}\left( x_{1}+\eta
_{0}\right) ^{-4} \\
&&-2\left( x_{1}+\eta _{0}\right) ^{2}(\beta _{0}-1)(n-1)+\delta \alpha
_{1}\left( x_{1}+\eta _{0}\right) m),
\end{eqnarray*}%
and the terms $D_{1}\left( \chi \varphi \right) ,$ $D_{2}\left( \chi \varphi
\right) $ are given by divergence forms which include the function $\varphi $
and are given in the proof explicitly.
\end{lemma}

The proof of Lemma 2 is technical and lenghty, and we postpone it to
Appendix.

In (8), the signs of the terms of $\left\vert \nabla _{^{\prime }x}\varphi
\right\vert ^{2}$ and $\left\vert \nabla _{y}\varphi \right\vert ^{2}$ are
different. Thus we need to perform another estimation:

\begin{lemma}
The following equality holds: 
\begin{eqnarray}
-(x_{1}+\eta _{0})\varphi (L_{0}\varphi )\chi ^{2} &=&\varphi
_{x_{1}}^{2}\chi ^{2}+\chi ^{2}\left( x_{1}+\eta _{0}\right) ^{2}\left(
|\nabla _{^{\prime }x}\varphi |^{2}-\sum_{i,j=1}^{m}a_{ij}\varphi
_{y_{i}}\varphi _{y_{j}}\right)  \notag \\
&&+\sigma _{2}(\lambda ,\delta )\varphi ^{2}\chi ^{2}+D_{3}\left( \varphi
\right)  \TCItag{9}
\end{eqnarray}%
for any function $\varphi \in C^{2}(\overline{\Omega })$. Here 
\begin{equation*}
\sigma _{2}(\lambda ,\delta )=\lambda ^{2}\nu ^{2}\sigma _{21}+\lambda \nu
\sigma _{22}+\sigma _{23},
\end{equation*}%
\begin{eqnarray*}
\sigma _{21} &=&-2\psi ^{-2\nu -2}\left( \delta ^{2}+\left( x_{1}+\eta
_{0}\right) ^{2}\left( \left\vert \nabla _{^{\prime }x}\psi \right\vert
^{2}-\sum_{i,j=1}^{m}a_{ij}\psi _{y_{i}}\psi _{y_{j}}\right) \right) , \\
\sigma _{22} &=&-(\nu +1)\psi ^{-\nu -2}\left( \delta ^{2}+\left( x_{1}+\eta
_{0}\right) ^{2}\left( \left\vert \nabla _{^{\prime }x}\psi \right\vert
^{2}-\sum_{i,j=1}^{m}a_{ij}\psi _{y_{i}}\psi _{y_{j}}\right) \right) \\
&&+\psi ^{-\nu -1}\left( x_{1}+\eta _{0}\right) ^{2}\left(
n-1-2\sum_{i,j=1}^{m}\frac{\partial a_{ij}}{\partial y_{j}}\psi
_{y_{i}}-\sum_{i=1}^{m}a_{ij}\right) , \\
\sigma _{23} &=&\frac{1}{2}\sum_{i,j=1}^{m}\frac{\partial ^{2}a_{ij}}{%
\partial y_{i}\partial y_{j}}(x+\eta _{0})^{2}
\end{eqnarray*}%
and%
\begin{eqnarray*}
D_{3}\left( \varphi \right) &=&-((\varphi \varphi _{x_{1}}+\lambda \nu
\delta \psi ^{-\nu -1}\varphi ^{2})\chi ^{2})_{x_{1}} \\
&&-\sum_{i=2}^{n}\frac{\partial }{\partial x_{i}}((\varphi \varphi
_{x_{i}}+\lambda \nu \psi _{x_{i}}\psi ^{-\nu -1}\varphi ^{2})\chi
^{2}\left( x_{1}+\eta _{0}\right) ^{2}) \\
&&+\sum_{i=1}^{m}\frac{\partial }{\partial y_{i}}\left( \sum_{j=1}^{m}\left(
a_{ij}\left( \varphi _{y_{j}}\varphi +\lambda \nu \psi _{y_{j}}\psi ^{-\nu
-1}\varphi ^{2}\right) -\frac{1}{2}\frac{\partial a_{ij}}{\partial y_{j}}%
\varphi ^{2}\right) \chi ^{2}\left( x_{1}+\eta _{0}\right) ^{2}\right) .
\end{eqnarray*}%
Lemma 3 can be proved easily by direct calculations and we omit the proof
here.
\end{lemma}

Now we will proceed to the completion of the proof of Lemma 1. We multiply
equality (9) by $2\lambda \nu \beta _{0}$ and add to inequality (8) to have%
\begin{eqnarray}
&&\psi ^{\nu +1}\left( L_{0}\varphi \right) ^{2}\chi ^{2}-2\lambda \nu \beta
_{0}(x_{1}+\eta _{0})\varphi \left( L_{0}\varphi \right) \chi ^{2}  \notag \\
&\geq &2\lambda \nu \delta \left( x_{1}+\eta _{0}\right) ^{-3}\varphi
_{x_{1}}^{2}\chi ^{2}+2\lambda \nu \left( x_{1}+\eta _{0}\right)
^{2}\left\vert \nabla _{^{\prime }x}\varphi \right\vert ^{2}\chi ^{2}  \notag
\\
&&+\lambda \nu \delta \alpha _{1}\left( x_{1}+\eta _{0}\right) \left\vert
\nabla _{y}\varphi \right\vert ^{2}\chi ^{2}-2\lambda \nu \beta _{0}\left(
x_{1}+\eta _{0}\right) ^{2}\chi ^{2}\sum_{i,j=1}^{m}a_{ij}\varphi
_{y_{i}}\varphi _{y_{j}}  \notag \\
&&+2\lambda ^{3}\nu ^{4}\delta ^{4}\left( x_{1}+\eta _{0}\right) ^{-2}\psi
^{-2\nu -3}\varphi ^{2}\chi ^{2}+\sigma _{3}(\lambda ,\delta )\varphi
^{2}\chi ^{2}  \notag \\
&&+D_{1}(\chi \varphi )+D_{2}(\chi \varphi )+2\lambda \nu \beta
_{0}D_{3}\left( \varphi \right) \text{,}  \TCItag{10}
\end{eqnarray}%
for $\delta >\delta _{0}$, $\lambda >\lambda _{0},$ where $\sigma
_{3}(\lambda ,\delta )=\sigma _{1}(\lambda ,\delta )+2\lambda \nu \beta
_{0}\sigma _{2}(\lambda ,\delta ).$ We set $\delta _{1}=\frac{2}{\alpha _{1}}%
\left( 1+\frac{3}{4}m\beta _{0}\gamma M\right) $.

Since%
\begin{equation*}
-\sum_{i,j=1}^{m}a_{ij}\varphi _{y_{i}}\varphi _{y_{j}}\geq -Mm\left\vert
\nabla _{y}\varphi \right\vert ^{2},
\end{equation*}%
we can estimate the coefficient of $\left\vert \nabla _{y}\varphi
\right\vert ^{2}$:%
\begin{eqnarray}
&&\lambda \nu \delta \alpha _{1}\left( x_{1}+\eta _{0}\right) \left\vert
\nabla _{y}\varphi \right\vert ^{2}\chi ^{2}-2\lambda \nu \beta _{0}\left(
x_{1}+\eta _{0}\right) ^{2}\sum_{i,j=1}^{m}a_{ij}\varphi _{y_{i}}\varphi
_{y_{j}}\chi ^{2}  \notag \\
&\geq &\lambda \nu \left( x_{1}+\eta _{0}\right) (\delta \alpha _{1}-2\beta
_{0}mM\left( x_{1}+\eta _{0}\right) )\left\vert \nabla _{y}\varphi
\right\vert ^{2}\chi ^{2}  \notag \\
&\geq &2\lambda \nu \left( x_{1}+\eta _{0}\right) \left\vert \nabla
_{y}\varphi \right\vert ^{2}\chi ^{2}  \TCItag{11}
\end{eqnarray}%
for $\delta \geq \delta _{1}.$

As for the coefficient of $\varphi ^{2}$, we can write $\sigma _{3}(\lambda
,\delta )$ in the form%
\begin{equation*}
\sigma _{3}(\lambda ,\delta )=\lambda ^{3}\nu ^{3}\sigma _{31}+\lambda
^{2}\nu ^{2}\sigma _{32}+\lambda \nu \sigma _{33},
\end{equation*}%
where $\sigma _{31}=\sigma _{11}+2\beta _{0}\sigma _{21},$ $\sigma
_{32}=\sigma _{12}+2\beta _{0}\sigma _{22},$ $\sigma _{33}=2\beta _{0}\sigma
_{23}.$

Since the functions $a_{ij},$ $\psi ,$ $\psi _{x_{i}},$ $\psi _{y_{j}}$ are
bounded in the space $C\left( \overline{\Omega }_{\gamma }\right) ,$ it is
clear that the function$\ \tilde{\sigma}_{31}=\frac{\sigma _{31}}{\delta
^{3}\nu \psi ^{-2\nu -3}}$ is bounded uniformly with respect to $\left(
x,y\right) \in \overline{\Omega }_{\gamma }:$ $\left\vert \tilde{\sigma}%
_{31}\right\vert \leq M_{1}$, $M_{1}>0.$ Then we see that%
\begin{eqnarray*}
\lambda ^{3}\nu ^{4}\delta ^{4}\left( x_{1}+\eta _{0}\right) ^{-2}\psi
^{-2\nu -3}+\lambda ^{3}\nu ^{3}\sigma _{31} &>&\lambda ^{3}\nu ^{4}\delta
^{4}\psi ^{-2\nu -3}\left( 1+\frac{1}{\delta }\tilde{\sigma}_{31}\right) \\
&\geq &\lambda ^{3}\nu ^{4}\delta ^{4}\psi ^{-2\nu -3}\left( 1-\frac{1}{%
\delta }M_{1}\right) \\
&=&\lambda ^{3}\nu ^{4}\delta ^{4}\psi ^{-2\nu -3}\left( 1-\frac{1}{\delta }%
\frac{\delta _{2}}{2}\right) \\
&\geq &\frac{1}{2}\lambda ^{3}\nu ^{4}\delta ^{4}\psi ^{-2\nu -3}
\end{eqnarray*}%
for $\delta \geq \delta _{2}=2M_{1}$. Here we note that $\left( x_{1}+\eta
_{0}\right) ^{-2}\geq (\frac{3}{4}\gamma )^{-2}>1.$

On the other hand, it is obvious that, for fixed $\delta \geq \delta _{2}$, $%
\nu >1$, the functions $\sigma _{32}$ and $\sigma _{33}$ are also bounded on 
$\Omega _{\gamma }$, that is, there exist constants $M_{2},$ $M_{3}>0$ such
that $\left\vert \sigma _{32}\right\vert \leq M_{2},$ $\left\vert \sigma
_{33}\right\vert \leq M_{3}.$ Thus we have%
\begin{eqnarray}
&&\lambda ^{3}\nu ^{4}\delta ^{4}\left( x_{1}+\eta _{0}\right) ^{-2}\psi
^{-2\nu -3}+\lambda ^{3}\nu ^{3}\sigma _{31}+\lambda ^{2}\nu ^{2}\sigma
_{32}+\lambda \nu \sigma _{33}  \notag \\
&\geq &\frac{1}{2}\lambda ^{3}\nu ^{4}\delta ^{4}\psi ^{-2\nu -3}-\lambda
^{2}\nu ^{2}M_{2}-\lambda \nu M_{3}  \notag \\
&\geq &0  \TCItag{12}
\end{eqnarray}%
for $\lambda \geq \lambda _{1}=\max \left\{ M_{2},\sqrt{M_{3}}\right\} ,$
which yields%
\begin{equation}
2\lambda ^{3}\nu ^{4}\delta ^{4}\left( x_{1}+\eta _{0}\right) ^{-2}\psi
^{-2\nu -3}-\sigma _{3}(\lambda ,\delta )\geq \lambda ^{3}\nu ^{4}\delta
^{4}\psi ^{-2\nu -3}.  \tag{13}
\end{equation}%
Consequently, inequalities (10), (11) and (13) imply that 
\begin{eqnarray*}
&&\psi ^{\nu +1}\left( L_{0}\varphi \right) ^{2}\chi ^{2}-2\lambda \nu \beta
_{0}\varphi \left( L_{0}\varphi \right) \chi ^{2} \\
&\geq &2\lambda \nu \delta \left( x_{1}+\eta _{0}\right) ^{-3}\varphi
_{1}^{2}\chi ^{2}+2\lambda \nu \left( x_{1}+\eta _{0}\right) ^{2}\chi
^{2}\left\vert \nabla _{^{\prime }x}\varphi \right\vert ^{2}\chi ^{2} \\
&&+2\lambda \nu \left( x_{1}+\eta _{0}\right) \left\vert \nabla _{y}\varphi
\right\vert ^{2}\chi ^{2}+\lambda ^{3}\nu ^{4}\delta ^{4}\psi ^{-2\nu
-3}\varphi ^{2}\chi ^{2}+D\left( \varphi \right) ,
\end{eqnarray*}%
for $\delta \geq \delta _{\ast }=\max \left\{ \delta _{0},\delta _{1},\delta
_{2}\right\} $ and $\lambda \geq \lambda _{\ast }=\max \left\{ \lambda
_{0},\lambda _{1}\right\} $, where $D\left( \varphi \right) =D_{1}\left(
\chi \varphi \right) +D_{2}\left( \chi \varphi \right) +2\lambda \nu \beta
_{0}D_{3}\left( \varphi \right) $. Thus the proof of Lemma 1 is complete.

\section{The Proof of Theorem 1}

Let $\left( u,g\right) $ be a solution to (1) - (3) with $u_{0}\equiv 0$ in $%
\Omega _{\gamma }$. Since $f\left( x,y^{\prime },0\right) \neq 0$ and $f\in
C^{2}\left( \overline{\Omega }\right) $, there exists a number $0<\gamma <1$
such that $f(x,y)\neq 0$ also in $\Omega _{\gamma }.$ We assume that $\gamma
,$ which was introduced before, satisfies this condition. We define a new
unknown function $w=\dfrac{u}{f}$ in $\Omega _{\gamma }.$ Then dividing
equation (5) by $f\left( x,y\right) $ and taking into account relations
(2)-(3), we obtain%
\begin{eqnarray}
&&\left( x_{1}+\eta _{0}\right) ^{-1}w_{x_{1}x_{1}}+\left( x_{1}+\eta
_{0}\right) (\Delta _{^{\prime }x}w-\overset{m}{\underset{i,j=1}{\sum }}%
a_{ij}w_{y_{i}y_{j}})  \notag \\
&&+\left( x_{1}+\eta _{0}\right) (\overset{n}{\underset{i=1}{\sum }}\bar{a}%
_{i}w_{x_{i}}+\overset{m}{\underset{j=1}{\sum }}\bar{b}_{j}w_{y_{j}}+\bar{a}%
_{0}w)  \notag \\
&=&\left( x_{1}+\eta _{0}\right) g,  \TCItag{14}
\end{eqnarray}%
\begin{equation}
w\left( 0,^{\prime }x,y\right) =w_{x_{1}}\left( 0,^{\prime }x,y\right) =0, 
\tag{15}
\end{equation}%
\begin{equation}
w\left( x,y^{\prime },0\right) =0,  \tag{16}
\end{equation}%
where%
\begin{eqnarray*}
\bar{a}_{0} &=&(\left( x_{1}+\eta _{0}\right) ^{-2}f_{x_{1}x_{1}}+\Delta
_{^{\prime
}x}f-\sum_{i,j=1}^{m}a_{ij}f_{y_{i}y_{j}}+\sum_{i=1}^{n}a_{i}f_{x_{i}} \\
&&+\sum_{j=1}^{m}b_{j}f_{y_{j}}+a_{0}f)f^{-1},
\end{eqnarray*}%
\begin{equation*}
\bar{a}_{1}=(2\left( x_{1}+\eta _{0}\right) ^{-2}f_{x_{1}}+a_{1}f)f^{-1},
\end{equation*}%
\begin{eqnarray*}
\bar{a}_{i} &=&(2f_{x_{i}}+a_{i}f)f^{-1},\text{ }i=2,...,n, \\
\bar{b}_{j}
&=&(-\sum_{i=1}^{m}(a_{ji}f_{y_{j}}+a_{ij}f_{y_{i}})+b_{j}f)f^{-1},\text{ }%
j=1,...,m.
\end{eqnarray*}%
Differentiating equation (14) with respect to $y_{m},$ setting $z=w_{y_{m}}$
and using (16), we obtain the integro-differential equation%
\begin{eqnarray}
&&\left( x_{1}+\eta _{0}\right) ^{-1}z_{x_{1}x_{1}}+(x_{1}+\eta _{0})\left(
\Delta _{^{\prime }x}z-\sum_{i,j=1}^{m}a_{ij}z_{y_{i}y_{j}}\right)  \notag \\
&&+(x_{1}+\eta _{0})\left( \sum\limits_{i=1}^{n}\bar{a}_{i}z_{x_{i}}+%
\sum_{j=1}^{m}\bar{b}_{j}z_{y_{j}}+\bar{a}_{0}z\right)  \notag \\
&&+(x_{1}+\eta _{0})\left( \sum_{i=1}^{n}\frac{\partial \bar{a}_{i}}{%
\partial y_{m}}I_{x_{i}}z+\sum_{j=1}^{m}\frac{\partial \bar{b}_{j}}{\partial
y_{m}}I_{y_{j}}z+\frac{\partial \bar{a}_{0}}{\partial y_{m}}Iz\right)  \notag
\\
&=&0  \TCItag{17}
\end{eqnarray}%
with the Cauchy data%
\begin{equation}
z\left( 0,^{\prime }x,y\right) =z_{x_{1}}\left( 0,^{\prime }x,y\right) =0, 
\tag{18}
\end{equation}%
where%
\begin{equation*}
Iz=\left\{ 
\begin{array}{c}
\int_{0}^{y_{m}}z(x,y^{\prime },\tau )d\tau ,\text{ }y_{m}\geq 0 \\ 
\int_{y_{m}}^{0}z(x,y^{\prime },\tau )d\tau ,\text{ }y_{m}<0%
\end{array}%
\right. ,
\end{equation*}%
\begin{equation*}
I_{x_{i}}z=\left\{ 
\begin{array}{c}
\frac{\partial }{\partial x_{i}}\int_{0}^{y_{m}}z(x,y^{\prime },\tau )d\tau ,%
\text{ }y_{m}\geq 0 \\ 
\frac{\partial }{\partial x_{i}}\int_{y_{m}}^{0}z(x,y^{\prime },\tau )d\tau ,%
\text{ }y_{m}<0%
\end{array}%
\right. ,\text{ }i=1,..,n;
\end{equation*}%
\begin{equation*}
I_{y_{j}}z=\left\{ 
\begin{array}{c}
\frac{\partial }{\partial y_{j}}\int_{0}^{y_{m}}z(x,y^{\prime },\tau )d\tau ,%
\text{ }y_{m}\geq 0 \\ 
\frac{\partial }{\partial y_{j}}\int_{y_{m}}^{0}z(x,y^{\prime },\tau )d\tau ,%
\text{ }y_{m}<0%
\end{array}%
\right. ,\text{ }j=1,...,m.
\end{equation*}%
We now prove that, if $z\left( x,y\right) $ satisfies (17) and (18), then $%
z\left( x,y\right) =0$ in $\Omega _{\gamma }$.

From (17), we obtain 
\begin{eqnarray}
\left( L_{0}z\right) ^{2} &=&(x_{1}+\eta _{0})^{2}\left(
\sum\limits_{i=1}^{n}\left( \bar{a}_{i}z_{x_{i}}+\frac{\partial \bar{a}_{i}}{%
\partial y_{m}}I_{x_{i}}z\right) +\sum_{j=1}^{m}\left( \bar{b}_{j}z_{y_{j}}+%
\frac{\partial \bar{b}_{j}}{\partial y_{m}}I_{y_{j}}z\right) \right.  \notag
\\
&&+\left. \bar{a}_{0}z+\frac{\partial \bar{a}_{0}}{\partial y_{m}}Iz\right)
^{2}  \notag \\
&\leq &6M_{5}\max \left\{ n,m\right\} \left( x_{1}+\eta _{0}\right)
^{2}\left( \left\vert \nabla _{x}z\right\vert
^{2}+\sum\limits_{i=1}^{n}(I_{x_{i}}z)^{2}+|\nabla _{y}z|^{2}\right.  \notag
\\
&&+\left. \sum_{j=1}^{m}(I_{y_{j}}z)^{2}+z^{2}+(Iz)^{2}\right) ,  \TCItag{19}
\end{eqnarray}%
where $M_{5}>0$ depends on $M$ and $\Vert f\Vert _{C^{2}(\overline{\Omega }%
_{\gamma })}$.

On the other hand, by Lemma 1, we can write%
\begin{eqnarray}
&&\left( L_{0}\varphi \right) ^{2}\chi ^{2}+(x_{1}+\eta
_{0})^{2}(L_{0}\varphi )^{2}\chi ^{2}+\lambda ^{2}\nu ^{2}\beta
_{0}^{2}\varphi ^{2}\chi ^{2}  \notag \\
&\geq &\psi ^{\nu +1}\left( L_{0}\varphi \right) ^{2}\chi ^{2}-2\lambda \nu
\beta _{0}\varphi (x_{1}+\eta _{0})(L_{0}\varphi )\chi ^{2}  \notag \\
&\geq &2\lambda \nu \delta \left( x_{1}+\eta _{0}\right) ^{-3}\varphi
_{x_{1}}^{2}\chi ^{2}+2\lambda \nu \left( x_{1}+\eta _{0}\right)
^{2}\left\vert \nabla _{^{\prime }x}\varphi \right\vert ^{2}\chi ^{2}  \notag
\\
&&+2\lambda \nu \left( x_{1}+\eta _{0}\right) \left\vert \nabla _{y}\varphi
\right\vert ^{2}\chi ^{2}+\lambda ^{3}\nu ^{4}\delta ^{4}(x_{1}+\eta
_{0})^{-2}\psi ^{-2\nu -3}\varphi ^{2}\chi ^{2}+D\left( \varphi \right) 
\TCItag{20}
\end{eqnarray}%
for $\delta >\delta _{\ast },$ $\lambda >\lambda _{\ast }$. Taking $\varphi
\equiv z$ in (20) and using (19), we obtain 
\begin{eqnarray}
&&6M_{5}\max \left\{ n,m\right\} \left( x_{1}+\eta _{0}\right)
^{2}(\left\vert \nabla _{x}z\right\vert
^{2}+\sum\limits_{i=1}^{n}(I_{x_{i}}z)^{2}  \notag \\
&&+\left\vert \nabla _{y}z\right\vert
^{2}+\sum_{j=1}^{m}(I_{y_{j}}z)^{2}+(Iz)^{2})\chi ^{2}(1+\left( x_{1}+\eta
_{0}\right) ^{2})  \notag \\
&&+(6M_{5}\max \left\{ n,m\right\} \left( x_{1}+\eta _{0}\right)
^{2}(1+\left( x_{1}+\eta _{0}\right) ^{2})+\lambda ^{2}\nu ^{2}\beta
_{0}^{2})z^{2}\chi ^{2}  \notag \\
&\geq &2\lambda \nu \delta \left( x_{1}+\eta _{0}\right)
^{-3}z_{x_{1}}^{2}\chi ^{2}+2\lambda \nu \left( x_{1}+\eta _{0}\right)
^{2}\left\vert \nabla _{^{\prime }x}z\right\vert ^{2}\chi ^{2}  \notag \\
&&+2\lambda \nu \left( x_{1}+\eta _{0}\right) \left\vert \nabla
_{y}z\right\vert ^{2}\chi ^{2}+\lambda ^{3}\nu ^{4}\delta ^{4}(x_{1}+\eta
_{0})^{-2}\psi ^{-2\nu -3}z^{2}\chi ^{2}+D\left( z\right)  \TCItag{21}
\end{eqnarray}%
Here we shall use the following lemma, whose proof is given in Appendix.

\begin{lemma}
The following relations hold:%
\begin{eqnarray*}
\int_{\Omega _{\gamma }}\left( Iz\right) ^{2}\chi ^{2}d\Omega _{\gamma }
&\leq &\gamma \int_{\Omega _{\gamma }}z^{2}\chi ^{2}d\Omega _{\gamma }, \\
\int_{\Omega _{\gamma }}\left( I_{x_{i}}z\right) ^{2}\chi ^{2}d\Omega
_{\gamma } &\leq &\gamma \int_{\Omega _{\gamma }}z_{x_{i}}^{2}\chi
^{2}d\Omega _{\gamma }, \\
\int_{\Omega _{\gamma }}\left( I_{y_{j}}z\right) ^{2}\chi ^{2}d\Omega
_{\gamma } &\leq &\gamma \int_{\Omega _{\gamma }}z_{y_{j}}^{2}\chi
^{2}d\Omega _{\gamma },
\end{eqnarray*}%
where $i=1,...,n;$ $j=1,...,m.$
\end{lemma}

Integrating inequality (21) on $\Omega _{\gamma }$ and using Lemma 4, we have%
\begin{eqnarray*}
&&6M_{5}\max \left\{ n,m\right\} (1+\gamma )\int_{\Omega _{\gamma
}}(1+\left( x_{1}+\eta _{0}\right) ^{2})\left( x_{1}+\eta _{0}\right) ^{2} \\
&&\times (\left\vert \nabla _{x}z\right\vert ^{2}+\left\vert \nabla
_{y}z\right\vert ^{2}+z^{2})\chi ^{2}d\Omega _{\gamma }+\lambda ^{2}\nu
^{2}\beta _{0}^{2}\int_{\Omega _{\gamma }}z^{2}\chi ^{2}d\Omega _{\gamma } \\
&\geq &2\lambda \nu \delta \left( x_{1}+\eta _{0}\right) ^{-3}\int_{\Omega
_{\gamma }}z_{x_{1}}^{2}\chi ^{2}d\Omega _{\gamma }+2\lambda \nu
\int_{\Omega _{\gamma }}\left( x_{1}+\eta _{0}\right) ^{2}\left\vert \nabla
_{^{\prime }x}z\right\vert ^{2}\chi ^{2}d\Omega _{\gamma } \\
&&+2\lambda \nu \int_{\Omega _{\gamma }}\left( x_{1}+\eta _{0}\right)
\left\vert \nabla _{y}z\right\vert ^{2}\chi ^{2}d\Omega _{\gamma } \\
&&+\lambda ^{3}\nu ^{4}\delta ^{4}\int_{\Omega _{\gamma }}\psi ^{-2\nu
-3}z^{2}\chi ^{2}d\Omega _{\gamma }+D\left( z\right) .
\end{eqnarray*}%
Hence, if $\lambda \geq \lambda _{\ast }=12M_{5}\max \left\{ n,m\right\}
(1+\gamma )>1$ and $\nu \geq \delta ^{-4}(1+\beta _{0}^{2}+(\frac{3}{4}%
\gamma )^{2})$, then we obtain%
\begin{eqnarray}
&&\int_{\Omega _{\gamma }}(\lambda ^{3}\nu ^{3}z^{2}+\lambda \nu
(z_{x_{1}}^{2}+\left( x_{1}+\eta _{0}\right) ^{2}\left\vert \nabla
_{^{\prime }x}z\right\vert ^{2}+\left( x_{1}+\eta _{0}\right) \left\vert
\nabla _{y}z\right\vert ^{2}))\chi ^{2}d\Omega _{\gamma }  \notag \\
&\leq &-\int_{\Omega _{\gamma }}D\left( z\right) d\Omega _{\gamma }. 
\TCItag{22}
\end{eqnarray}%
Passing to the limit as $\lambda \rightarrow \infty $ in (22), we conclude
that%
\begin{equation*}
\int_{\Omega _{\gamma }}z^{2}d\Omega _{\gamma }\leq 0,
\end{equation*}%
which means that $z=0$ in $\Omega _{\gamma }$.

Varying the point\ $x_{0}=\left(
0,x_{2}^{0},x_{3}^{0},...,x_{n+m}^{0}\right) $ of the plane $x_{1}=0$, we
establish that $z=0$\ on $\tilde{\Omega}_{\gamma }=\left\{ (x,y)\in \Omega ;%
\text{ }0\leq \delta x_{1}\leq \gamma \right\} ,$\ that is $\frac{\partial w%
}{\partial y_{m}}=0$ on $\tilde{\Omega}_{\gamma }$. Then from equation (14)
by condition (16) we conclude that $g\left( x,y^{\prime }\right) =0$ on\ $%
\tilde{\Omega}_{\gamma }^{\prime }=\left\{ (x,y^{\prime })\in D\times
G^{\prime };\text{ }0\leq \delta x_{1}\leq \gamma \right\} $, \ where $%
\tilde{\Omega}_{\gamma }^{\prime }$\ the is projection of $\Omega $ onto $%
\mathbb{R}
^{n+m-1}$.$\ $Repeating the same argument, we see that $z=0$\ in $\tilde{%
\Omega}_{2\gamma }$ and $g\left( x,y^{\prime }\right) =0$\ on $\tilde{\Omega}%
_{2\gamma }^{\prime }$. Thus, continuing the argument, we complete the proof.

\section{Appendix}

\subsection{Proof of Lemma 2}

We introduce a new function%
\begin{equation*}
\vartheta =\chi \varphi .
\end{equation*}%
Using the relations%
\begin{equation*}
\varphi _{x_{1}x_{1}}=\chi ^{-1}\left( \vartheta _{x_{1}x_{1}}+2\lambda \nu
\psi ^{-\nu -1}\psi _{x_{1}}\vartheta _{x_{1}}+\delta ^{2}\phi _{1}\vartheta
\right) 
\end{equation*}%
\begin{equation*}
\Delta _{^{\prime }x}\varphi =\chi ^{-1}\left( \Delta _{^{\prime
}x}\vartheta +2\lambda \nu \psi ^{-\nu -1}(\nabla _{^{\prime }x}\psi ,\nabla
_{^{\prime }x}\vartheta )+\phi _{2}\vartheta \right) ,
\end{equation*}%
\begin{equation*}
\sum_{i,j=1}^{m}a_{ij}\varphi _{y_{i}y_{j}}=\chi ^{-1}\left(
\sum_{i,j=1}^{m}a_{ij}\left( \vartheta _{y_{i}y_{j}}+\lambda \nu \psi ^{-\nu
-1}\psi _{y_{i}}\vartheta _{y_{j}}+\lambda \nu \psi ^{-\nu -1}\psi
_{y_{j}}\vartheta _{y_{i}}\right) +\phi _{3}\vartheta \right) ,
\end{equation*}%
we obtain%
\begin{eqnarray}
&&\psi ^{\nu +1}\left( L_{0}\varphi \right) ^{2}\chi ^{2}  \notag \\
&=&\psi ^{\nu +1}\biggl\{\left( x_{1}+\eta _{0}\right) ^{-1}\vartheta
_{x_{1}x_{1}}+\left( x_{1}+\eta _{0}\right) \left( \Delta _{^{\prime
}x}\vartheta -\overset{m}{\underset{i,j=1}{\sum }}a_{ij}\vartheta
_{y_{i}y_{j}}\right)   \notag \\
&&+\vartheta \left( \left( x_{1}+\eta _{0}\right) ^{-1}\delta ^{2}\phi
_{1}+\left( x_{1}+\eta _{0}\right) \left( \phi _{2}-\phi _{3}\right)
)\right.   \notag \\
&&+2\lambda \nu \psi ^{-\nu -1}\left( \delta \left( x_{1}+\eta _{0}\right)
^{-1}\vartheta _{x_{1}}+\left( x_{1}+\eta _{0}\right) \left( (\nabla
_{^{\prime }x}\psi ,\nabla _{^{\prime }x}\vartheta )-\overset{m}{\underset{%
i,j=1}{\sum }}a_{ij}\psi _{y_{i}}\vartheta _{y_{j}}\right) \right) \biggr\}%
^{2}  \notag \\
&\geq &4\lambda \nu \biggl\{\left( x_{1}+\eta _{0}\right) ^{-1}\vartheta
_{x_{1}x_{1}}+\left( x_{1}+\eta _{0}\right) \left( \Delta _{^{\prime
}x}\vartheta -\overset{m}{\underset{i,j=1}{\sum }}a_{ij}\vartheta
_{y_{i}y_{j}}\right)   \notag \\
&&+\vartheta (\left( x_{1}+\eta _{0}\right) ^{-1}\delta ^{2}\phi _{1}+\left(
x_{1}+\eta _{0}\right) \left( \phi _{2}-\phi _{3}\right) )\biggr\}  \notag \\
&&\times \left( \delta \left( x_{1}+\eta _{0}\right) ^{-1}\vartheta
_{x_{1}}+\left( x_{1}+\eta _{0}\right) \left( (\nabla _{^{\prime }x}\psi
,\nabla _{^{\prime }x}\vartheta )-\overset{m}{\underset{i,j=1}{\sum }}%
a_{ij}\psi _{y_{i}}\vartheta _{y_{j}}\right) \right)   \notag \\
&=&:\sum\limits_{k=1}^{14}T_{k}\text{,}  \TCItag{23}
\end{eqnarray}%
where we set%
\begin{eqnarray*}
\phi _{1} &:&=\phi _{1}\left( \lambda ,\nu ,\psi \right) =\lambda ^{2}\nu
^{2}\psi ^{-2\nu -2}-\lambda \nu (\nu +1)\psi ^{-\nu -2}, \\
\phi _{2} &:&=\phi _{2}\left( \lambda ,\nu ,\psi \right) =\lambda ^{2}\nu
^{2}\phi _{21}\left( \nu ,\psi \right) -\lambda \nu \phi _{22}\left( \nu
,\psi \right) ,\text{ } \\
\phi _{3} &:&=\phi _{3}\left( \lambda ,\nu ,\psi \right) =\lambda ^{2}\nu
^{2}\phi _{31}\left( \nu ,\psi \right) -\lambda \nu \phi _{32}\left( \nu
,\psi \right) 
\end{eqnarray*}%
and 
\begin{eqnarray*}
\phi _{21}\left( \nu ,\psi \right)  &:&=\psi ^{-2\nu -2}\left\vert \nabla
_{^{\prime }x}\psi \right\vert ^{2},\text{ }\phi _{22}\left( \nu ,\psi
\right) =(\nu +1)\psi ^{-\nu -2}\left\vert \nabla _{^{\prime }x}\psi
\right\vert ^{2}-\left( n-1\right) \psi ^{-\nu -1} \\
\phi _{31}\left( \nu ,\psi \right)  &:&=\psi ^{-2\nu
-2}\sum_{i,j=1}^{m}a_{ij}\psi _{y_{i}}\psi _{y_{j}}, \\
\phi _{32}\left( \nu ,\psi \right)  &:&=\sum_{i,j=1}^{m}a_{ij}\left( ((\nu
+1)\psi ^{-\nu -2}-\psi ^{-\nu -1})\psi _{y_{i}}\psi _{y_{j}}\right) 
\end{eqnarray*}%
Noting that $\left\Vert a_{ij}\right\Vert _{C^{1}\left( \Omega \right) }\leq
M$, $\left\vert \psi _{x_{i}}\right\vert \leq \sqrt{2\gamma },\ \left( 2\leq
i\leq n\right) $ and $|\psi _{y_{k}}|\leq \sqrt{2\gamma }$\ $\left( 1\leq
k\leq m\right) $\ in$\ \Omega _{\gamma },$ we estimate the terms $T_{i},$ $%
1\leq i\leq 14$ as follows:%
\begin{equation}
T_{1}=4\lambda \nu \delta \left( x_{1}+\eta _{0}\right) ^{-2}\vartheta
_{x_{1}}\vartheta _{x_{1}x_{1}}=d_{1}\left( \vartheta \right) +4\lambda \nu
\delta \left( x_{1}+\eta _{0}\right) ^{-3}\vartheta _{x_{1}}^{2},  \tag{24}
\end{equation}%
where $d_{1}\left( \vartheta \right) =2\lambda \nu \delta (\left( x_{1}+\eta
_{0}\right) ^{-2}\vartheta _{x_{1}}^{2})_{x_{1}}$;%
\begin{eqnarray}
T_{2} &=&4\lambda \nu \delta \vartheta _{x_{1}}\Delta _{^{\prime
}x}\vartheta   \notag \\
&=&4\lambda \nu \delta \sum_{i=2}^{n}(\left( \vartheta _{x_{1}}\vartheta
_{x_{i}}\right) _{x_{i}}-2\lambda \nu \delta \left( \vartheta
_{x_{i}}^{2}\right) _{x_{1}})  \notag \\
&=&:d_{2}\left( \vartheta \right) \text{;}  \TCItag{25}
\end{eqnarray}%
\begin{eqnarray}
T_{3} &=&-4\lambda \nu \delta \vartheta _{x_{1}}\overset{m}{\underset{i,j=1}{%
\sum }}a_{ij}\vartheta _{y_{i}y_{j}}  \notag \\
&=&d_{3}\left( \vartheta \right) -2\lambda \nu \delta \overset{m}{\underset{%
i,j=1}{\sum }}(-2\frac{\partial a_{ij}}{\partial y_{j}}\vartheta
_{y_{i}}\vartheta _{x_{1}}+\frac{\partial a_{ij}}{\partial x_{1}}\vartheta
_{y_{i}}\vartheta _{y_{j}})  \notag \\
&\geq &d_{3}\left( \vartheta \right) -4\lambda \nu \delta \overset{m}{%
\underset{i,j=1}{\sum }}\left\vert \frac{\partial a_{ij}}{\partial y_{j}}%
\right\vert \left\vert \vartheta _{y_{i}}\vartheta _{x_{1}}\right\vert
+2\lambda \nu \delta \alpha _{1}\left\vert \nabla _{y}\vartheta \right\vert
^{2}  \notag \\
&\geq &d_{3}\left( \vartheta \right) +2\lambda \nu \delta (\alpha
_{1}-m(x_{1}+\eta _{0})\varepsilon _{0})\left\vert \nabla _{y}\vartheta
\right\vert ^{2}-2\lambda \nu \delta \frac{m^{2}M^{2}}{\varepsilon
_{0}\left( x_{1}+\eta _{0}\right) }\vartheta _{x_{1}}^{2},  \TCItag{26}
\end{eqnarray}%
where $d_{3}\left( \vartheta \right) =-2\lambda \nu \delta \overset{m}{%
\underset{i,j=1}{\sum }}(2\left( a_{ij}\vartheta _{y_{i}}\vartheta
_{x_{1}}\right) _{y_{j}}-\left( a_{ij}\vartheta _{y_{i}}\vartheta
_{y_{j}}\right) _{x_{1}}).$

Next 
\begin{eqnarray}
T_{4} &=&4\lambda \nu \delta ^{3}\left( x_{1}+\eta _{0}\right)
^{-2}\vartheta \vartheta _{x_{1}}\phi _{1}  \notag \\
&=&d_{4}\left( \vartheta \right) +4\lambda \nu \delta ^{3}\vartheta
^{2}\left( x_{1}+\eta _{0}\right) ^{-3}\phi _{1}  \notag \\
&&+4\lambda \nu \delta ^{4}\vartheta ^{2}(\nu +1)\left( x_{1}+\eta
_{0}\right) ^{-2}\phi _{4},  \TCItag{27}
\end{eqnarray}%
where $d_{4}\left( \vartheta \right) =2\lambda \nu \delta ^{3}(\left(
x_{1}+\eta _{0}\right) ^{-2}\vartheta ^{2}\phi _{1})_{x_{1}}$ and we set%
\begin{equation*}
\phi _{4}:=\phi _{4}\left( \lambda ,\nu ,\psi \right) =\lambda ^{2}\nu
^{2}\psi ^{-2\nu -3}-\frac{1}{2}\lambda \nu (\nu +2)\psi ^{-\nu -3}.
\end{equation*}%
We have 
\begin{eqnarray}
T_{5} &=&4\lambda \nu \delta \vartheta _{x_{1}}\vartheta \phi _{x_{1}} 
\notag \\
&=&d_{5}\left( \vartheta \right) +4\lambda \nu \delta ^{2}(\nu +1)\vartheta
^{2}\left( \phi _{4}\sum_{i=2}^{n}\psi _{i}^{2}+\frac{1}{2}\frac{(n-1)}{(\nu
+1)}\lambda \nu \psi ^{-\nu -2}\right) ,  \TCItag{28}
\end{eqnarray}%
where $d_{5}\left( \vartheta \right) =2\lambda \nu \delta \left( \vartheta
^{2}\phi _{x_{1}}\right) _{x_{1}};$%
\begin{equation}
T_{6}=-4\lambda \nu \delta \vartheta _{x_{1}}\vartheta \phi _{3}=d_{6}\left(
\vartheta \right) +2\lambda \nu \delta \vartheta ^{2}\left( \phi _{3}\right)
_{x_{1}},  \tag{29}
\end{equation}%
where $d_{6}\left( \vartheta \right) =-2\lambda \nu \delta \left( \vartheta
^{2}\phi _{3}\right) _{x_{1}};$%
\begin{equation}
T_{7}=4\lambda \nu \sum_{i=2}^{n}\psi _{x_{i}}\vartheta _{x_{i}}\vartheta
_{x_{1}x_{1}}=d_{7}\left( \vartheta \right) +2\lambda \nu \vartheta
_{x_{1}}^{2},  \tag{30}
\end{equation}%
where $d_{7}\left( \vartheta \right) =4\lambda \nu
\sum\limits_{i=2}^{n}(\left( \psi _{x_{i}}\vartheta _{x_{i}}\vartheta
_{x_{1}}\right) _{x_{1}}-2\lambda \nu \left( \psi _{x_{i}}\vartheta
_{x_{1}}^{2}\right) _{x_{i}});$%
\begin{eqnarray}
T_{8} &=&4\lambda \nu \left( x_{1}+\eta _{0}\right) ^{2}\sum_{i,j=2}^{n}\psi
_{x_{i}}\vartheta _{x_{i}}\vartheta _{x_{j}x_{j}}  \notag \\
&=&d_{8}\left( \vartheta \right) -4\lambda \nu \left( x_{1}+\eta _{0}\right)
^{2}\sum_{i,j=2}^{n}\delta _{ij}\vartheta _{x_{i}}\vartheta
_{x_{j}}+2\lambda \nu \left( x_{1}+\eta _{0}\right)
^{2}\sum_{i,j=2}^{n}\vartheta _{x_{j}}^{2}  \notag \\
&=&d_{8}\left( \vartheta \right) -2\lambda \nu \left( x_{1}+\eta _{0}\right)
^{2}\sum_{i=2}^{n}\vartheta _{x_{i}}^{2}(2-(n-1)),  \TCItag{31}
\end{eqnarray}%
where $d_{8}(\vartheta )=2\lambda \nu \sum\limits_{i,j=2}^{n}\left( 2(\psi
_{x_{i}}\vartheta _{x_{i}}\vartheta _{x_{j}}\left( x_{1}+\eta _{0}\right)
^{2})_{x_{j}}-\left( x_{1}+\eta _{0}\right) ^{2}\left( \psi
_{x_{i}}\vartheta _{x_{j}}^{2}\right) _{x_{i}}\right) .$ 
\begin{eqnarray}
T_{9} &=&-4\lambda \nu \left( x_{1}+\eta _{0}\right) ^{2}\sum_{i=2}^{n}\psi
_{x_{i}}\vartheta _{x_{i}}\sum_{k,s=1}^{m}a_{ks}\vartheta _{y_{k}y_{s}} 
\notag \\
&=&d_{9}\left( \vartheta \right) +2\lambda \nu \left( x_{1}+\eta _{0}\right)
^{2}\sum_{i=2}^{n}\sum_{k,s=1}^{m}(2\left( \psi _{x_{i}}a_{ks}\right)
_{y_{s}}\vartheta _{y_{k}}\vartheta _{x_{i}}-\left( \psi
_{x_{i}}a_{ks}\right) _{x_{i}}\vartheta _{y_{k}}\vartheta _{y_{s}})  \notag
\\
&\geq &-2\lambda \nu \left( x_{1}+\eta _{0}\right)
^{2}\sum_{i=2}^{n}\sum_{k,s=1}^{m}\left( 2\left\vert \left( \psi
_{x_{i}}a_{ks}\right) _{y_{s}}\vartheta _{y_{k}}\vartheta
_{x_{i}}\right\vert +\left\vert \left( \psi _{x_{i}}a_{ks}\right)
_{x_{i}}\vartheta _{y_{k}}\vartheta _{y_{s}}\right\vert \right)   \notag \\
&\geq &-2\lambda \nu \left( x_{1}+\eta _{0}\right) ^{2}M(1+\sqrt{2\gamma }%
)\left( 2mn\left\vert \nabla _{y}\vartheta \right\vert ^{2}+m^{2}\left\vert
\nabla _{^{\prime }x}\vartheta \right\vert ^{2}\right) ,  \TCItag{32}
\end{eqnarray}%
where $d_{9}\left( \vartheta \right) =-2\lambda \nu \left( x_{1}+\eta
_{0}\right) ^{2}\sum\limits_{k,s=1}^{m}\sum\limits_{k,s=1}^{m}(2\left( \psi
_{x_{i}}a_{ks}\vartheta _{x_{i}}\vartheta _{y_{k}}\right) _{s}-\left( \psi
_{x_{i}}a_{ks}\vartheta _{y_{k}}\vartheta _{y_{s}}\right) _{x_{i}});$%
\begin{eqnarray}
T_{10} &=&4\lambda \nu \sum_{i=2}^{n}\psi _{x_{i}}\vartheta
_{x_{i}}\vartheta \left( x_{1}+\eta _{0}\right) (\left( x_{1}+\eta
_{0}\right) ^{-1}\delta ^{2}\phi _{1}+\left( x_{1}+\eta _{0}\right) \left(
\phi _{2}-\phi _{3}\right) )  \notag \\
&=&d_{10}\left( \vartheta \right) -2\lambda \nu \vartheta
^{2}\sum_{i=2}^{n}(\psi _{x_{i}}\delta ^{2}\phi _{1}+\left( x_{1}+\eta
_{0}\right) ^{2}\psi _{x_{i}}(\phi _{2}-\phi _{3}))\text{, }  \TCItag{33}
\end{eqnarray}%
where $d_{10}\left( \vartheta \right) =2\lambda \nu
\sum\limits_{i=2}^{n}(\vartheta ^{2}\psi _{x_{i}}\left( x_{1}+\eta
_{0}\right) (\left( x_{1}+\eta _{0}\right) ^{-1}\delta ^{2}\phi _{1}+\left(
x_{1}+\eta _{0}\right) (\phi _{2}-\phi _{3})))_{x_{i}};$%
\begin{eqnarray}
T_{11} &=&-4\lambda \nu \sum_{i,j=1}^{m}a_{ij}\psi _{y_{i}}\vartheta
_{y_{j}}\vartheta _{x_{1}x_{1}}=d_{11}\left( \vartheta \right)   \notag \\
&&+4\lambda \nu \sum_{i,j=1}^{m}\frac{\partial a_{ij}}{\partial x_{1}}\psi
_{y_{i}}\vartheta _{y_{j}}\vartheta _{x_{1}}-2\lambda \nu
\sum_{i,j=1}^{m}\left( \frac{\partial a_{ij}}{\partial y_{j}}\psi
_{y_{i}}+a_{ij}\delta _{ij}\right) \vartheta _{x_{1}}^{2}  \notag \\
&\geq &-4\lambda \nu \sum_{i,j=1}^{m}\left\vert \frac{\partial a_{ij}}{%
\partial x_{1}}\psi _{y_{i}}\vartheta _{y_{j}}\vartheta _{x_{1}}\right\vert
-2\lambda \nu \sum_{i,j=1}^{m}\left\vert \left( \frac{\partial a_{ij}}{%
\partial y_{j}}\psi _{y_{i}}+a_{ij}\delta _{ij}\right) \right\vert \vartheta
_{x_{1}}^{2}  \notag \\
&\geq &-2\lambda \nu \left( x_{1}+\eta _{0}\right) \sqrt{2\gamma }%
m\left\vert \nabla _{y}\vartheta \right\vert ^{2}-2\lambda \nu M\left(
x_{1}+\eta _{0}\right) ^{-1}m^{2}\upsilon _{x_{1}}^{2}  \notag \\
&&-2\lambda \nu M\sqrt{2\gamma }m^{2}\vartheta _{x_{1}}^{2}-2\lambda \nu
Mm\vartheta _{x_{1}}^{2},  \TCItag{34}
\end{eqnarray}%
where $d_{11}\left( \vartheta \right) =-4\lambda \nu
\sum\limits_{i,j=1}^{m}(\left( a_{ij}\psi _{y_{i}}\vartheta
_{y_{j}}\vartheta _{x_{1}}\right) _{x_{1}}+2\lambda \nu \left( a_{ij}\psi
_{y_{i}}\vartheta _{x_{1}}^{2}\right) _{y_{j}});$%
\begin{eqnarray}
T_{12} &=&-4\lambda \nu \left( x_{1}+\eta _{0}\right)
^{2}\sum\limits_{i,j=1}^{m}a_{ij}\psi _{y_{i}}\vartheta
_{y_{j}}\sum\limits_{s=2}^{n}\vartheta _{x_{s}x_{s}}  \notag \\
&=&d_{12}\left( \vartheta \right) +4\lambda \nu \left( x_{1}+\eta
_{0}\right) ^{2}\sum\limits_{i,j=1}^{m}\sum\limits_{s=2}^{n}\frac{\partial
a_{ij}}{\partial x_{s}}\psi _{y_{i}}\vartheta _{y_{j}}\vartheta _{x_{s}} 
\notag \\
&&-2\lambda \nu \left( x_{1}+\eta _{0}\right)
^{2}\sum\limits_{i,j=1}^{m}\sum\limits_{s=2}^{n}\left( \frac{\partial a_{ij}%
}{\partial y_{j}}\psi _{y_{i}}+a_{ij}\delta _{ij}\right) \vartheta
_{x_{s}}^{2}  \notag \\
&\geq &-2\lambda \nu \left( x_{1}+\eta _{0}\right) ^{2}M\sqrt{2\gamma }%
\left( m^{2}\left\vert \nabla _{^{\prime }x}\vartheta \right\vert
^{2}+(n-1)m\left\vert \nabla _{y}\vartheta \right\vert ^{2}\right)   \notag
\\
&&-2\lambda \nu \left( x_{1}+\eta _{0}\right) ^{2}Mm(\sqrt{2\gamma }%
m+1)\left\vert \nabla _{^{\prime }x}\vartheta \right\vert ^{2},  \TCItag{35}
\end{eqnarray}%
where $d_{12}\left( \vartheta \right) =-4\lambda \nu (\left( x_{1}+\eta
_{0}\right) ^{2}\sum\limits_{i,j=1}^{m}\sum\limits_{s=2}^{n}\left(
a_{ij}\psi _{y_{i}}\vartheta _{y_{j}}\vartheta _{x_{s}}\right)
_{x_{s}}+2\lambda \nu \left( x_{1}+\eta _{0}\right) ^{2}\left( a_{ij}\psi
_{y_{i}}\vartheta _{x_{s}}^{2}\right) _{y_{j}});$%
\begin{eqnarray}
T_{13} &=&4\lambda \nu \left( x_{1}+\eta _{0}\right)
^{2}\sum\limits_{i,j,k,s=1}^{m}a_{ij}\psi _{y_{i}}\vartheta
_{y_{j}}a_{ks}\vartheta _{y_{k}y_{s}}=d_{13}\left( \vartheta \right)   \notag
\\
&&-2\lambda \nu \left( x_{1}+\eta _{0}\right)
^{2}\sum_{i,j,k,s=1}^{m}(2\left( a_{ij}a_{ks}\psi _{y_{i}}\right)
_{y_{s}}\vartheta _{y_{j}}\vartheta _{y_{k}}-\left( a_{ij}a_{ks}\psi
_{y_{i}}\right) _{y_{j}}\vartheta _{y_{k}}\vartheta _{y_{s}})  \notag \\
&\geq &-2\lambda \nu \left( x_{1}+\eta _{0}\right) ^{2}\overset{m}{\underset{%
i,j,k,s=1}{\sum }}(2\left\vert \left( a_{ij}a_{\kappa s}\psi _{y_{i}}\right)
_{y_{s}}\right\vert \left\vert \vartheta _{y_{j}}\vartheta
_{y_{k}}\right\vert +\left\vert \left( a_{ij}a_{\kappa s}\psi
_{y_{i}}\right) _{y_{j}}\right\vert \left\vert \vartheta _{y_{k}}\vartheta
_{y_{s}}\right\vert )  \notag \\
&\geq &-6\lambda \nu \left( x_{1}+\eta _{0}\right) ^{2}M^{2}\left( 2\sqrt{%
2\gamma }+1\right) m^{3}\left\vert \nabla _{y}\vartheta \right\vert ^{2}, 
\TCItag{36}
\end{eqnarray}%
where $d_{13}\left( \vartheta \right) =2\lambda \nu \left( x_{1}+\eta
_{0}\right) ^{2}\overset{m}{\underset{i,j,k,s=1}{\sum }}(2(a_{ij}\psi
_{y_{i}}\vartheta _{y_{j}}\vartheta _{y_{k}}a_{ks})_{y_{s}}-\left(
a_{ij}a_{ks}\psi _{y_{i}}\vartheta _{y_{k}}\vartheta _{y_{s}}\right)
_{y_{j}});$%
\begin{eqnarray}
T_{14} &=&-4\lambda \nu \left( x_{1}+\eta _{0}\right)
\sum_{i,j=1}^{m}a_{ij}\psi _{y_{i}}\vartheta _{y_{j}}\vartheta (\left(
x_{1}+\eta _{0}\right) ^{-1}\delta ^{2}\phi _{1}+\left( x_{1}+\eta
_{0}\right) (\phi _{2}-\phi _{3}))  \notag \\
&=&d_{14}\left( \vartheta \right) +2\lambda \nu \vartheta
^{2}\sum_{i,j=1}^{m}(a_{ij}\psi _{y_{i}}(\delta ^{2}\phi _{1}+\left(
x_{1}+\eta _{0}\right) ^{2}(\phi _{2}-\phi _{3})))_{y_{j}},  \TCItag{37}
\end{eqnarray}%
where $d_{14}\left( \vartheta \right) =-2\lambda \nu
\sum\limits_{i,j=1}^{m}(a_{ij}\psi _{y_{i}}\vartheta ^{2}\left( x_{1}+\eta
_{0}\right) (\left( x_{1}+\eta _{0}\right) ^{-1}\delta ^{2}\phi _{1}+\left(
x_{1}+\eta _{0}\right) (\phi _{2}-\phi _{3})))_{y_{j}}.$ Then by relations
(24)-(37), we see that%
\begin{eqnarray}
&&\psi ^{\nu +1}\left( L_{0}\varphi \right) ^{2}\chi ^{2}\geq 2\lambda \nu
\delta \left( x_{1}+\eta _{0}\right) ^{-3}\beta _{1}\vartheta
_{1}^{2}-2\lambda \nu \left( x_{1}+\eta _{0}\right) ^{2}(\beta
_{0}-1)\left\vert \nabla _{^{\prime }x}\vartheta \right\vert ^{2}  \notag \\
+ &&2\lambda \nu \left( x_{1}+\eta _{0}\right) \beta _{2}\left\vert \nabla
_{y}\vartheta \right\vert ^{2}+\left( \lambda ^{3}\nu ^{3}\beta _{3}+\lambda
^{2}\nu ^{2}\beta _{4}\right) \vartheta ^{2}+D_{1}\left( \vartheta \right) ,
\TCItag{38}
\end{eqnarray}%
where%
\begin{eqnarray*}
\beta _{1} &=&2-Mm\left( x_{1}+\eta _{0}\right) ^{2}(Mm\varepsilon
_{0}^{-1}-m\delta ^{-1}-\delta ^{-1}(\sqrt{2\gamma }m+1)(x_{1}+\eta _{0})),
\\
\beta _{2} &=&\delta (\alpha _{1}-m\varepsilon _{0})-\beta _{21}, \\
\beta _{21} &=&\sqrt{2\gamma }m+(2n(1+\sqrt{2\gamma })+(n-1)\sqrt{2\gamma }%
+3m^{2}M(2\sqrt{2\gamma }+1))\left( x_{1}+\eta _{0}\right) mM
\end{eqnarray*}%
\begin{eqnarray*}
\beta _{3} &=&4\delta ^{3}\left( x_{1}+\eta _{0}\right) ^{-3}\psi ^{-2\nu
-2}+4\delta ^{4}(\nu +1)\left( x_{1}+\eta _{0}\right) ^{-2}\psi ^{-2\nu -3}
\\
&&+4\delta ^{2}(\nu +1)\psi ^{-2\nu -3}\left\vert \nabla _{^{\prime }x}\psi
\right\vert ^{2}+\beta _{31}, \\
\beta _{31} &=&2\delta (\phi _{31})_{x_{1}}+2\left\vert \nabla _{^{\prime
}x}\psi \right\vert ^{2}(\left( x_{1}+\eta _{0}\right) ^{2}(\phi _{31}-\phi
_{21})-\delta ^{2}\psi ^{-2\nu -2}\phi _{1}) \\
&&+2\sum_{i,j=1}^{m}(a_{ij}\psi _{y_{i}}(\left( x_{1}+\eta _{0}\right)
^{2}(\phi _{21}-\phi _{31})+\delta ^{2}\psi ^{-2\nu -2}))_{y_{j}}, \\
\beta _{4} &=&-4\delta ^{3}\left( x_{1}+\eta _{0}\right) ^{-3}(\nu +1)\psi
^{-\nu -2}+4\delta ^{2}(\nu +1)(\frac{1}{2}\frac{(n-1)}{(\nu +1)}\psi ^{-\nu
-2} \\
&&-\frac{1}{2}(\nu +2)\psi ^{-\nu -3}\left\vert \nabla _{^{\prime }x}\psi
\right\vert ^{2}-\delta ^{2}\left( x_{1}+\eta _{0}\right) ^{-2}\frac{1}{2}%
(\nu +2)\psi ^{-\nu -3})-2\delta (\phi _{32})_{x_{1}} \\
&&+2\left\vert \nabla _{^{\prime }x}\psi \right\vert ^{2}(\left( x_{1}+\eta
_{0}\right) ^{2}(\phi _{22}-\phi _{32})+\delta ^{2}(\nu +1)\psi ^{-\nu
-2}\phi _{1}) \\
&&+2\sum_{i,j=1}^{m}(a_{ij}\psi _{y_{i}}(\left( x_{1}+\eta _{0}\right)
^{2}(\phi _{32}-\phi _{22})-\delta ^{2}(\nu +1)\psi ^{-\nu -2}))_{y_{j}}
\end{eqnarray*}%
and $D_{1}\left( \vartheta \right) =\sum\limits_{k=1}^{14}d_{k}\left(
\vartheta \right) .$

Now we shall evaluate the expressions $\beta _{1},$ $\beta _{2},$ $\beta
_{3},$ $\beta _{4}$ in (38), respectively.

Since $\delta \geq 4,$ $x_{1}+\eta _{0}<\frac{3}{4}\gamma ,$ and $Mm\left( 
\frac{3}{4}\gamma \right) ^{2}(Mm\varepsilon _{0}^{-1}+m+(\sqrt{2\gamma }m+1)%
\frac{3}{4}\gamma )<1$ by (6), we obtain%
\begin{eqnarray*}
\beta _{1} &>&2-Mm\left( x_{1}+\eta _{0}\right) ^{2}(Mm\varepsilon
_{0}^{-1}+m+(\sqrt{2\gamma }m+1)(x_{1}+\eta _{0})) \\
&>&2-Mm\left( \frac{3}{4}\gamma \right) ^{2}\left( Mm\varepsilon
_{0}^{-1}+m+(\sqrt{2\gamma }m+1)\frac{3}{4}\gamma \right) >1,
\end{eqnarray*}%
which implies%
\begin{eqnarray}
&&2\lambda \nu \delta \left( x_{1}+\eta _{0}\right) ^{-3}\vartheta
_{x_{1}}^{2}(2-M^{2}m^{2}\varepsilon _{0}^{-1}(x_{1}+\eta
_{0})^{2}-Mm^{2}\delta ^{-1}(x_{1}+\eta _{0})^{2}  \notag \\
&&-M\delta ^{-1}m(\sqrt{2\gamma }m+1)\left( x_{1}+\eta _{0}\right) ^{3}) 
\notag \\
&\geq &2\lambda \nu \delta \left( x_{1}+\eta _{0}\right) ^{-3}\vartheta
_{x_{1}}^{2}\text{.}  \TCItag{39}
\end{eqnarray}%
Next, by choosing $0<\varepsilon _{0}<\frac{\alpha _{1}}{4m},$ and setting $%
\delta _{3}=\frac{4m}{\alpha _{1}}\sqrt{2\gamma }l_{1},$ we see that%
\begin{equation}
\beta _{2}>\delta \left( \alpha _{1}-m\varepsilon _{0}\right) -m\sqrt{%
2\gamma }l_{1}\geq \frac{1}{2}\delta \alpha _{1}  \tag{40}
\end{equation}%
for $\delta \geq \delta _{3},$ where 
\begin{equation*}
l_{1}=1+(n-1)M\frac{3}{4}\gamma +(2nM+3M^{2}m^{2})(1+2\sqrt{2\gamma })\frac{3%
}{4}\left( \frac{\gamma }{2}\right) ^{1/2}.
\end{equation*}%
It is clear that%
\begin{eqnarray}
\beta _{3} &>&4\delta ^{4}(\nu +1)\left( x_{1}+\eta _{0}\right) ^{-2}\psi
^{-2\nu -3}-\left\vert \beta _{31}\right\vert   \notag \\
&=&4\delta ^{4}(\nu +1)\left( x_{1}+\eta _{0}\right) ^{-2}\psi ^{-2\nu -3}(1-%
\frac{1}{\delta ^{2}}\left\vert \tilde{\beta}_{31}\right\vert )  \TCItag{41}
\end{eqnarray}%
where $\tilde{\beta}_{31}=\frac{1}{4\delta ^{2}}\left( x_{1}+\eta
_{0}\right) ^{2}\psi ^{2\nu +3}\frac{1}{\nu +1}\beta _{31}.$

Since the functions $a_{ij}$ and $\psi $ bounded in the space $C^{2}\left( 
\overline{\Omega }_{\gamma }\right) ,$ the function $\tilde{\beta}_{31}$ is
bounded uniformly with respect to $x\in \overline{\Omega }_{\gamma },$ that
is there exist a number $M_{4}>0$ such that $\left\vert \tilde{\beta}%
_{31}\right\vert \leq M_{4}.$ Then%
\begin{equation}
\beta _{3}>2\delta ^{4}(\nu +1)\left( x_{1}+\eta _{0}\right) ^{-2}\psi
^{-2\nu -3}  \tag{42}
\end{equation}%
for $\delta \geq \delta _{4}=\sqrt{2M_{4}}.$

By the same reasons, $\beta _{4}$ is$\ $also bounded on $\overline{\Omega }%
_{\gamma }$:$\ $that is $\left\vert \beta _{4}\right\vert \leq M_{5}$, $%
M_{5}>0.$ Moreover, since%
\begin{equation*}
2\lambda ^{3}\nu ^{3}-\lambda ^{2}\nu ^{2}M_{5}\geq 0
\end{equation*}%
for $\lambda \geq \lambda _{2}=M_{5}\left( \delta ,\nu ,\gamma \right) ,$ $%
\nu >1,$ we can write%
\begin{eqnarray}
&&\lambda ^{3}\nu ^{3}\beta _{3}+\lambda ^{2}\nu ^{2}\beta _{4}\geq 2\lambda
^{3}\nu ^{3}\delta ^{4}(\nu +1)\left( x_{1}+\eta _{0}\right) ^{-2}\psi
^{-2\nu -3}-\lambda ^{2}\nu ^{2}M_{2}  \notag \\
&=&2\lambda ^{3}\nu ^{4}\delta ^{4}\left( x_{1}+\eta _{0}\right) ^{-2}\psi
^{-2\nu -3}+2\lambda ^{3}\nu ^{3}\delta ^{4}\left( x_{1}+\eta _{0}\right)
^{-2}\psi ^{-2\nu -3}-\lambda ^{2}\nu ^{2}M_{2}  \notag \\
&\geq &2\lambda ^{3}\nu ^{4}\delta ^{4}\left( x_{1}+\eta _{0}\right)
^{-2}\psi ^{-2\nu -3}  \TCItag{43}
\end{eqnarray}%
for $\delta \geq \delta _{0}=\max \left\{ 4,\delta _{3},\delta _{4}\right\} ,
$ $\nu >1,$ $\lambda \geq \lambda _{1}.$

Thus, by inequalities (39)-(43), we have%
\begin{eqnarray}
&&\psi ^{\nu +1}\left( L_{0}\varphi \right) ^{2}\chi ^{2}\geq 2\lambda \nu
\delta \left( x_{1}+\eta _{0}\right) ^{-3}\vartheta _{x_{1}}^{2}-2\lambda
\nu \left( x_{1}+\eta _{0}\right) ^{2}(\beta _{0}-1)\left\vert \nabla
_{^{\prime }x}\vartheta \right\vert ^{2}  \notag \\
&&+\lambda \nu \delta \alpha _{1}\left( x_{1}+\eta _{0}\right) \left\vert
\nabla _{y}\vartheta \right\vert ^{2}+2\lambda ^{3}\nu ^{4}\delta ^{4}\left(
x_{1}+\eta _{0}\right) ^{-2}\psi ^{-2\nu -3}\vartheta ^{2}+d\left( \vartheta
\right)  \TCItag{44}
\end{eqnarray}%
for $\delta \geq \delta _{0}$,$\ \lambda \geq \lambda _{1}$,$\ \nu >1.$

Finally, taking into account the equality $\vartheta =\chi \varphi $ and the
relations%
\begin{equation*}
\vartheta _{x_{1}}^{2}=\varphi _{x_{1}}^{2}\chi ^{2}-\lambda ^{2}\nu
^{2}\delta ^{2}\psi ^{-2\nu -2}\varphi ^{2}\chi ^{2}-\lambda \nu (\nu
+1)\delta ^{2}\psi ^{-\nu -2}\varphi ^{2}\chi ^{2}-\lambda \nu \left( \delta
\psi ^{-\nu -1}\varphi ^{2}\chi ^{2}\right) _{x_{1}},
\end{equation*}%
\begin{eqnarray*}
\left\vert \nabla _{^{\prime }x}\vartheta \right\vert ^{2} &=&\left\vert
\nabla _{^{\prime }x}\varphi \right\vert ^{2}\chi ^{2}-\lambda ^{2}\nu
^{2}\left\vert \nabla _{^{\prime }x}\psi \right\vert ^{2}\psi ^{-2\nu
-2}\varphi ^{2}\chi ^{2}-\lambda \nu \left( \nu +1\right) \left\vert \nabla
_{^{\prime }x}\psi \right\vert ^{2}\psi ^{-\nu -2}\varphi ^{2}\chi ^{2} \\
&&+\lambda \nu \left( n-1\right) \psi ^{-\nu -1}\varphi ^{2}\chi
^{2}-\lambda \nu \dsum\limits_{i=2}^{n}\left( \psi _{x_{i}}\psi ^{-\nu
-1}\varphi ^{2}\chi ^{2}\right) _{x_{i}}, \\
\left\vert \nabla _{y}\vartheta \right\vert ^{2} &=&\left\vert \nabla
_{y}\varphi \right\vert ^{2}\chi ^{2}-\lambda ^{2}\nu ^{2}\left\vert \nabla
_{y}\psi \right\vert ^{2}\psi ^{-2\nu -2}\varphi ^{2}\chi ^{2}-\lambda \nu
\left( \nu +1\right) \left\vert \nabla _{y}\psi \right\vert ^{2}\psi ^{-\nu
-2}\varphi ^{2}\chi ^{2} \\
&&+\lambda \nu m\psi ^{-\nu -1}\varphi ^{2}\chi ^{2}-\lambda \nu
\dsum\limits_{k=1}^{m}\left( \psi _{y_{k}}\psi ^{-\nu -1}\varphi ^{2}\chi
^{2}\right) _{y_{k}},
\end{eqnarray*}%
from (44) we obtain (8), where 
\begin{eqnarray}
D_{2}(\chi \varphi ) &=&-2\lambda ^{2}\nu ^{2}\delta ^{2}\left( \left(
x_{1}+\eta _{0}\right) ^{-3}\varphi ^{2}\psi ^{-\nu -1}\chi ^{2}\right)
_{x_{1}}  \notag \\
&&-\lambda ^{2}\nu ^{2}\delta \alpha _{1}\sum_{i=1}^{m}\left( \left(
x_{1}+\eta _{0}\right) \varphi ^{2}\psi _{y_{i}}\psi ^{-\nu -1}\chi
^{2}\right) _{y_{i}}  \notag \\
&&+2\lambda ^{2}\nu ^{2}\left( x_{1}+\eta _{0}\right) ^{2}(\beta
_{0}-1)\sum_{i=2}^{n}\left( \varphi ^{2}\psi _{i}\psi ^{-\nu -1}\chi
^{2}\right) _{x_{i}}.  \TCItag{45}
\end{eqnarray}%
Thus the proof of Lemma 2 is complete.

\subsection{Proof of Lemma 4}

Let $K$ be the part of the domain $\Omega _{\gamma }$ produced by the plane $%
y_{m}=0$ and $\sigma (x,y^{\prime })$ is the distance between the boundary
point $(x,y)$ of $\Omega _{\gamma }$ and the point $(x,y^{\prime },0)\in K.$

Since $\chi ^{2}$ is monotonic with respect to $y_{m}$, we have%
\begin{eqnarray*}
&&\int_{\Omega _{\gamma }}\chi ^{2}(Iz)^{2}d\Omega _{\gamma } \\
&=&\int_{K}\left(\int\nolimits_{0}^{\sigma \left( x,y^{\prime }\right) }\chi
^{2}(Iz)^{2}dy_{m}+\int\nolimits_{-\sigma \left( x,y^{\prime }\right)
}^{0}\chi ^{2}(Iz)^{2}dy_{m}\right)dxdy^{\prime } \\
&\leq &\int_{K}\left(\sqrt{\gamma }\int\nolimits_{0}^{a}%
\int_{0}^{y_{m}}z^{2}(x,y^{\prime },\tau )\chi ^{2}d\tau dy_{m}+\sqrt{\gamma 
}\int\nolimits_{-a}^{0}\int_{y_{m}}^{0}z^{2}(x,y^{\prime },\tau )\chi
^{2}d\tau dy_{m}\right)dxdy^{\prime } \\
&\leq &\gamma \int_{K}(\int\nolimits_{0}^{a}z^{2}(x,y^{\prime },\tau )\chi
^{2}(x,y^{\prime },\tau )d\tau +\int\nolimits_{-a}^{0}z^{2}(x,y^{\prime
},\tau )\chi ^{2}(x,y^{\prime },\tau )d\tau )dxdy^{\prime } \\
&=&\gamma \int_{\Omega _{\gamma }}z^{2}\chi ^{2}d\Omega _{\gamma }
\end{eqnarray*}%
for $a=\sigma (x,y^{\prime})$. Here we used the relations%
\begin{eqnarray*}
\int\nolimits_{0}^{a}\chi ^{2}(x,y)\int_{0}^{y_{m}}z^{2}(x,y^{\prime },\tau
)d\tau dy_{m} &=&\int\nolimits_{0}^{a}\int_{0}^{y_{m}}\chi
^{2}(x,y)z^{2}(x,y^{\prime },\tau )d\tau dy_{m} \\
&=&\int\nolimits_{0}^{a}z^{2}(x,y^{\prime },\tau )\int_{\tau }^{a}\chi
^{2}dy_{m}d\tau \\
&\leq &\int\nolimits_{0}^{a}z^{2}(x,y^{\prime },\tau )\chi ^{2}(x,y^{\prime
},\tau )(a-\tau )d\tau \\
&\leq &\sqrt{\gamma }\int\nolimits_{0}^{a}z^{2}(x,y^{\prime },\tau )\chi
^{2}(x,y^{\prime },\tau )d\tau
\end{eqnarray*}%
and 
\begin{eqnarray*}
\int\nolimits_{-a}^{0}\chi ^{2}\int_{y_{m}}^{0}z^{2}(x,y^{\prime },\tau
)d\tau dy_{m} &=&\int\nolimits_{-a}^{0}z^{2}(x,y^{\prime },\tau
)\int_{-a}^{\tau }\chi ^{2}dy_{m}d\tau \\
&\leq &\sqrt{\gamma }\int\nolimits_{-a}^{0}z^{2}(x,y^{\prime },\tau )\chi
^{2}(x,y^{\prime },\tau )d\tau .
\end{eqnarray*}%
Similarly we can prove the rest part. 

\end{document}